\documentclass[a4paper,12pt]{article}

\usepackage{amsmath,amssymb,theorem}

{\theorembodyfont{\rm}
  \newtheorem{defi}{Definition}[section]
  \newtheorem{rem}[defi]{Remark}

}
  
  \newtheorem{prop}[defi]{Proposition}
  \newtheorem{thm}[defi]{Theorem}

\newcommand\Nset{\mathbb N}

\begin{document}

\title{Some Constructions of Divisible Designs from Laguerre Geometries}

\author{Sabine Giese \and Hans Havlicek \and Ralph-Hardo Schulz}

\date{}

\maketitle

\begin{abstract}
In the nineties, A.G. Spera introduced a construction principle for
divisible  designs.  Using this method, we get series of divisible
designs from finite  Laguerre  geometries. We show a close connection
between some of these divisible designs  and  divisible designs whose
construction was based on a conic in a plane  of  a  3-dimensional
projective space.
\end{abstract}


\section{Introduction}

One interesting sort of designs is that of
 divisible  designs. In 1992, A.\ G.\ Spera introduced a method to
construct these designs by  $t-R$-homogeneous ($t-R$-transitive)
$R$-permutation groups  (\cite{sp1}).
Here, $R$ denotes an equivalence relation on the elements of a finite
set.
Once we have such an $R$-permutation group acting on a finite set $X$,
the main problem is the calculation of the parameters of the divisible
design. Especially the determination of the order of the stabilizer
of a chosen base block in the $R$-permutation group, which is needed for the
calculation of the parameter $\lambda$, is often not trivial.
For that purpose, we have to obtain suitable conditions on the construction.
There are already several known
examples of constructions using Spera's construction principle
(\cite{sp1}, \cite{schsp1}, \cite{cesch1},\cite{cesch2}). In  1999,
C. Cerroni and R.-H. Schulz gave one such construction starting from
a conic in a plane of the 3-dimensional projective space PG$(3,q)$
(\cite{cesch1}).\\

In this paper, we construct several series of new 3-divisible designs
again using
this method  but starting from a Laguerre geometry
$\Sigma(\mathrm{GF}(q),\mathbb{D}(\mathrm{GF}(q)))$ with
$\mathbb{D}(\mathrm{GF}(q))$ being the ring  of dual numbers over the
finite field GF($q$).  We show a close connection between some of these
designs and those constructed by Cerroni and Schulz. We will
see that the initial situations of both constructions are mutually
dual and that the parameters of our first series of 3-divisible
designs of Theorem \ref{meinDD} are equal to the parameters of the
series of 3-divisible designs by Cerroni and Schulz (cf. Theorem
\ref{schulz/cerroni}) mentioned above. The other series of Theorem
\ref{meinDD} seem to be new.\\

\section{Preliminary terms, definitions and results}

\subsection{About Divisible Designs}

Let $X$ be a finite set with an equivalence relation $R$ on its
elements. We  denote  by $[x]$ the $R$-equivalence class containing $x
\in X$ and define $S:=\{[x]| x \in X\}$.

\begin{defi}
A subset $Y$ of $X$ is called \emph{$R$-transversal}, if $|Y \cap[x]|
\le 1 \mbox{   for all }\\ x \in X$.
\end{defi}

\begin{defi}[Definition of divisible designs]
Let $t,s,k,\lambda_t$ be positive integers with $t \le k <v=|X|$.
A triple $\mathcal{D}=(X,\mathcal{B},S)$ is called
$t-(s,k,\lambda_t)$-\emph{divisible design} (or
$t-(s,k,\lambda_t)$-\textup{DD})  if\\

\begin{tabular}{ll}
(1) & $\mathcal{B}$ is a set of $R$-transversal subsets of $X$ with
    $|B|=k$ for  all  $B \in\mathcal{B}$;\\
(2) & $|[x]|=s$ for all $x \in X$;\\
(3) & for every $R$-transversal $t$-subset $Y$ of $X$ there exist
exactly   $\lambda_t$\\
& ele\-ments  of $\mathcal{B}$ containing $Y$.
\end{tabular}\\

\noindent The elements of $X$ are called $points$, those of
$\mathcal{B}\;  blocks$, and the elements of $S\; point\; classes$.
\end{defi}

In this paper, we always suppose that every divisible design is
\emph{simple},  that  means that there exist no repeated blocks. Note
that a $t$-divisible  design is also a
$(t-1)$-divisible design with $\lambda_{t- 1}=\lambda_t(v-st+s)(k-
t+1)^{-1}$.
We shall use this observation and consider as well the 2-divisible designs
arising from the 3-divisible designs to be constructed below.\\
One way of constructing divisible designs is given by the
following proposition of A.G. Spera.

\subsubsection{SPERA's construction principle}

\begin{defi}
Let $G$ be a group acting on the set $X$ and $R$ an equivalence
relation on $X$  which is $G$-invariant, that is,
\[xRy \Longrightarrow x^gRy^g (\mbox{ for all } g \in G, \;\; x,y\in X);\]
then $\Lambda=(G,X,R)$ is called an $R-$\emph{group}.
(The group $G$ induces a permutation group on $X$, but not necessarily
faithfully.)
\end{defi}

\begin{defi}
$\Lambda$ is called \emph{$t-R$-transitive} if for any two $R$-transversal
$t$-tuples $(x_1,\ldots,x_t)$ and $(y_1,\ldots,y_t)$ of elements of
$X$ there exists an element $g$ of $G$ such that $y_i=x_i^g \,\mbox{
  for }   i=1,2,\ldots,t$.
\end{defi}

\begin{prop}[A.G. Spera, \cite{sp1}] \label{spera}
Let $\Lambda=(G,X,R)$ be a finite $t-R$-transitive
$R$-group, and let $B$ be an $R$-transversal subset of $X$ with
$t \le k:=|B|<v:=|X|$,
then the incidence structure
$ D(\Lambda, B)=(X,B^G,S) \mbox{~ for ~} B^G=\{B^g \,|\, g\in G\}$
is a $t-(s,k,\lambda_t)$-divisible design with $s=|[x]|$ for some $x \in X$,
$k=|B|$,
\[b=\frac{|G|}{|G_B|} \mbox{ and } \lambda_t=|G|\left(\begin{array}{c} k\\t
\end{array}\right)\left(|G_B|\left(\begin{array}{c} vs^{-1} \\ t
  \end{array}\right)s^t\right)^{-1},\]
where $G_B$ denotes the setwise
stabiliser  of  $B$ and $b$ the number of blocks of $D(\Lambda, B)$.
Moreover, $G$ induces a point- and block-transitive automorphism group of
$D(\Lambda,B)$.
\end{prop}

\begin{rem}\label{conDD}
To construct divisible designs by using Spera's proposition, we need a
finite set   $X$ with an equivalence relation $R$ on its elements and a
finite $t$-$R$-transitive   $R$-permutation group acting on this
set. Then, we have to choose a so called   'base block' and
calculate the parameters.
\end{rem}

By using Spera's proposition, C. Cerroni and R.-H. Schulz constructed
the following series of divisible  designs \cite{cesch1}.

\begin{thm}[Cerroni, Schulz]\label{schulz/cerroni}
Let $q=p^n$, where $p$ is a prime, and
let $n,i\in\Nset$ with $i|n$.
If $q$ is odd, there exists a $3-(q,p^i+1,1)$ - \textup{DD} with
$q^2+q$ points, having as a point-  and block-transitive automorphism
group $T\tilde G$ with $\tilde G\cong$  \textup{GO}$(3,q)$ and $T$ the
translation group of \textup{AG}$(3,q)$.
\end{thm}

By starting from PGO$(3,q)$, which acts
3-transitively on a given conic in a plane of the 3-dimensional projective
space PG$(3,q)$, Cerroni and Schulz con\-struc\-ted a
3-$R$-transitive $R$-permutation group ($R$ denotes the parallelism
re\-lation) of a finite set of affine planes in the
corresponding\footnote{The 3-dimensional affine space whose ideal
  plane is the plane considered above containing the given conic.}
3-dimensional affine space.
After choosing a base block, they used Spera's proposition to construct these
divisible  designs.\\
The whole construction and the proof can be found in \cite{cesch1}.
In this paper, we will describe the idea of the construction in part 4 in
order to compare it with our construction below.\\

\subsection{A Laguerre Geometry over the finite field $\mathrm{GF}(q)$}

\begin{defi}
A \emph{dual number over the finite field GF$(q)$} is an ordered pair
$(a,b)$ with   $a,b \in \mathrm{GF}(q)$ and with the following properties.
\end{defi}

Two dual numbers are equal if their components are equal.
The rules for addition and multiplication are:

\begin{enumerate}
\item[(i)] $(a,b)+(a', b'):=(a+a', b+b')$

\item[(ii)] $(a,b)\cdot(a', b'):=(a a',a b' + b a').$
\end{enumerate}

A dual number $(a,b)$ can also be represented either by a matrix
$\bigl( \begin{smallmatrix} a&b \\ 0&a \end{smallmatrix}\bigr)$
 with $ a,b \in \mathrm{GF}(q)$
or in the following form:
$a+b\epsilon$  where $\epsilon$ is any chosen element satisfying
$\epsilon^2 = 0$, for instance $\epsilon=(0,1)$.
The rules for addition and multiplication correspond to those
for matrices.

\begin{rem}
As a subring of the ring of matrices, the set of dual numbers with the
given addition and multiplication is a ring with $1$. We denote this
ring by $\mathbb{D}(\mathrm{GF}(q))$.
\end{rem}

In analogy to Benz (\cite{benz1}, p.24), we have

\begin{prop}
$\mathbb{D}(\mathrm{GF}(q))$ is a commutative local ring.
\end{prop}

The only maximal ideal $N$ contains all non-invertible elements of
$\mathbb{D}(\mathrm{GF}(q))$. For any ring R, we define $\mbox{R}^*$
as the multiplicative group of all
invertible elements. Here
$\mathcal{R}:=\mathbb{D}(\mathrm{GF}(q))^*=\mathbb{D}(\mathrm{GF}(q))\setminus
N$.

\begin{samepage}
\begin{defi}[Laguerre algebra]
For $K$ a field, a $K$-algebra $A$ is called a \emph{Laguerre algebra}
provided  there exists a two-sided ideal $M$ of $A$ with $A^*=A\setminus
M$ and $A=K  \oplus  M$.
\end{defi}
\end{samepage}

$\mathbb{D}(\mathrm{GF}(q))$ is commutative, hence the above defined
ideal $N$ is two-sided.
Furthermore, $\mathbb{D}(\mathrm{GF}(q))$ is the
direct sum of the embedded field $\mathrm{GF}(q)$ with $N$, where
$\mathrm{GF}(q)$ is
identified with the set of diagonal matrices
$\bigl( \begin{smallmatrix} a&0 \\ 0&a \end{smallmatrix}\bigr)$
(or the set of pairs $(a,0)$ or $a+0\epsilon$ with $a\in
\mathrm{GF}(q)$, depending
on the manner of representation), hence we get the well known result

\begin{prop}\label{lagalg}
$\mathbb{D}(\mathrm{GF}(q))$ is a Laguerre algebra.
\end{prop}

\begin{defi}
We define the \emph{projective line
  $\mathbb{P}(\mathbb{D}(\mathrm{GF}(q)))$  over
  $\mathbb{D}(\mathrm{GF}(q))$} as the set of
all equivalence classes of admissible pairs. Here, we call a pair
$(x_1,x_2),\, x_1,x_2 \in \mathbb{D}(\mathrm{GF}(q))$ an
\emph{admissible pair over $\mathbb{D}(\mathrm{GF}(q))$}
if at least one element is invertible.
Two admissible pairs $(x_1,x_2),(y_1,y_2)$ are called
\emph{equivalent} if there   exists an element $r \in \mathcal{R}$
such that $x_i=ry_i, \, i=1,2$.  We  call the elements of
$\mathbb{P}(\mathbb{D}(\mathrm{GF}(q)))$ \emph{points}.
Since $\mathbb{D}(\mathrm{GF}(q))$ is a local ring, this definition of
the projective line over $\mathbb{D}(\mathrm{GF}(q))$ is equivalent to
that given by Herzer \cite{herz} on p. 785.
\end{defi}

\begin{defi}
Two points $P,Q \in \mathbb{P}(\mathbb{D}(\mathrm{GF}(q)))$ with
$P=\mathcal{R}(p_1,p_2),
Q=\mathcal{R}(q_1,q_2),\;p_i,q_i\in\mathbb{D}(\mathrm{GF}(q)),\,
i=1,2$ are called \emph{parallel} if
\[p_1q_2-q_1p_2= \left| \begin{array}{cc}   p_1&p_2 \\ q_1&q_2
  \end{array} \right| \notin \mathcal{R}.\]
\end{defi}

In analogy to the more general definition of chain geometry by Benz
in  (\cite{benz1}, p. 94), we define, in this paper, the chain
geometry $\Sigma  (K,D)$ with $K:=\mathrm{GF}(q),
D:=\mathbb{D}(\mathrm{GF}(q))$ as an incidence
structure whose points are the  elements of $\mathbb{P}(D)$ and whose blocks
(chains) are the images of $\mathbb{P}(K)$ under the  projective group of
$\mathbb{P}(D)$ (cf. Def. \ref{gamma} and  \cite{herz}, p. 790).
Since $D$ is a Laguerre algebra, the following holds:

\begin{prop} \textup{(\cite{herz})}
 $\Sigma(K, D)$ is a so called \emph{Laguerre geometry}, i.e., the
 parallelism relation is an equivalence relation on $\mathbb{P}(D)$ and every
 chain of $\Sigma(K, D)$ meets every parallel class of points.
\end{prop}

Now consider the chain geometry $\Sigma(K, D)$ whose points are the
elements  of the  projective line over $D$. We can partition the points
of $\mathbb{P}(D)$ into  \emph{proper}  and \emph{improper} points
depending on the invertiblity of the  second component  or,
equivalently, on the parallelism to the point $\mathcal{R}(1,0)$.
Every proper point can  be represented as $\mathcal{R}(p,1),\, p\in D$
and every improper  point as $\mathcal{R}(1,\delta\epsilon),\, \delta\in K$.

\begin{prop}
\begin{enumerate}
\item[(i)] Let $P,Q \in \mathbb{P}(D)$ be two proper points with
  $P=\mathcal{R}(p_1+p_2\epsilon,1)$ and
  $Q=\mathcal{R}(q_1+q_2\epsilon,1)$. They are parallel
  iff $p_1 = q_1$.

\item[(ii)] Improper points are always parallel.

\item[(iii)] A proper point is never parallel to an improper point.
\end{enumerate}
\end{prop}

The easy proof of this uses
\[P||Q \Leftrightarrow \bigl| \begin{smallmatrix} p_1+p_2\epsilon & 1
  \\ q_1+q_2\epsilon &1 \end{smallmatrix}\bigr|
= p_1-q_1+(p_2-q_2)\epsilon  \notin \mathcal{R}  \Leftrightarrow
p_1-q_1=0.\]

\begin{rem}
By the parallelism relation, we get $q+1$ equivalence classes with
$q$ elements each: \\
\centerline{
$~\{\mathcal{R}(x+b\epsilon, 1)| b \in K\}$, with $x\in K$ and
$[\mathcal{R}(1,0)]=\{\mathcal{R}(1, \delta\epsilon)| \delta \in K\}$.}
\end{rem}

\begin{rem}\label{emGF}
We can embed the projective line over $K$ into $\mathbb{P}(D)$. The
elements of $\mathbb{P}(K)$ form the following transversal subset of
$\mathbb{P}(D)$: \\
\centerline{
$\tilde
K:=\{\mathcal{R}(p_1+0\epsilon, 1)| p_1 \in  K\} \cup \mathcal{R}(1,0)$.}
\end{rem}

\begin{defi}\label{gamma}
One defines the projective group of $\mathbb{P}(D)$ as the group of
all regular $2\times 2$  matrices with entries in $\mathbb{D}(k)$
factorised by the subgroup \\
$\{\bigl( \begin{smallmatrix} r&0 \\ 0&r \end{smallmatrix}\bigr)
|\, r\in\mathcal{R}\}.$
We denote it by $\Gamma(D)$. (cf. \cite{benz1})
\end{defi}

\begin{prop} \textup{(\cite{benz1})} \label{transi}
$\Gamma(D)$ acts sharply 3-$R$-transitively on the point set of
$\mathbb{P}(D)$ and preserves parallelism.
\end{prop}

By Remark \ref{conDD}, we are now able to construct a divisible design
using Proposition \ref{spera}.\\

\section{Construction of divisible designs from a La\-guerre geometry}

\begin{thm}\label{meinDD}
Let $n,i\in\Nset$ with $i|n$ and let $q=p^n$, where $p$ is a prime.
Then, there exist $3-(q,k,\lambda_3)$-divisible designs, each
with $q(q+1)$ points with the parameters $k$ and $\lambda_3$ given
in Table 1, where $p$ and $i$ are subject to the conditions given
there. These $3$-divisible designs admit
$\Gamma(\mathbb{D}(\mathrm{GF}(q)))$ as a
point- and block-transitive automorphism group. The same holds for
the corresponding $2$-divisible designs.
\end{thm}

\begin{table}[h]
\begin{center}
\caption{}
\begin{tabular}{|c|c|c|l|}
  \hline
No.\footnotemark &  $k$ & $\lambda_3$ & Conditions \\
  \hline\hline
(i)&  $p^i+1$ &$1$&  \\
  \hline
(ii)&  $p^i$   & $p^i-2$& $p^i>2$\\
  \hline
(iii)&  $p^i-1$ & $\frac{1}{2}(p^i-2)(p^i-3))$ &$p^i >3$\\
  \hline
(iv)&  $p^i-2$ & $\frac{1}{6}(p^i-2)(p^i-3)(p^i-4))$ & $p^i>4$\\
  \hline
&  $p^i-3$ & $\frac{1}{4}(p^i-3)(p^i-4)(p^i-5))$ & \\
  \cline{2-3}
&  $4$ & $6$ &  \raisebox{1.5ex}[-1.5ex]{$p^i>7$}\\
  \cline{2-4}
&  $p^i-3$ & $\frac{1}{24}(p^i-3)(p^i-4)(p^i-5))$ &\\
  \cline{2-3}
(v)&  $4$ & $1$ &  \raisebox{1.5ex}[-1.5ex]{$p=3$ and $p^i>5$}\\
  \cline{2-4}
&  $p^i-3$ & $\frac{1}{8}(p^i-3)(p^i-4)(p^i-5))$ &\\
  \cline{2-3}
&  $4$ & $3$ &  \raisebox{1.5ex}[-1.5ex]{$p>3$ and $p^i >5$}\\
  \cline{2-4}
&  $p^i-3$ & $\frac{1}{12}(p^i-3)(p^i-4)(p^i-5))$ & \\
   \cline{2-3}
&  $4$ & $2$ & \raisebox{1.5ex}[-1.5ex]{$p^i \equiv 1\mod 3$ and $p^i >5$}\\
  \hline
\end{tabular}
\end{center}
\end{table}

\begin{paragraph}{\it Proof of Theorem \ref{meinDD}:}
Consider the chain geometry
$\Sigma(\mathrm{GF}(q),\mathbb{D}(\mathrm{GF}(q)))$, $q=p^n$ where $p$
is a prime. We use the same notation as above. By
Prop. \ref{transi}, we have, with   $\Gamma(D)$, a 3-$R$-transitive
$R$-permutation group acting on the point set which consists of
the $q^2$ proper points and $q$ improper points. They are  divided
into $q+1$ parallel classes with $q$ elements each, giving the points
and the point classes of a DD.
By using the sharp
$3-R$-transitivity and determining the order of the orbit of a
transversal triple, we obtain $|\Gamma(D)| = q^4(q^2-1)$.
Now we determine the order of the stabiliser of the considered
base block which we choose for the different cases.\\

\noindent Let $i\in \Nset$ with $i|n$ where $q=p^n$.

\begin{itemize}
\item[(i)]

For $L:=\mathrm{GF}(p^i)$, we embed the projective line over $L$
\[\mathbb{P}(L)=\textup{PG}(1,p^i)=:B\]
into the projective line over $D$ and define it as our base block.\\
Notice that a projectivity of $\Gamma(D)$, which maps three distinct
points of $B$ onto points of $B$, belongs to $\Gamma(L)$
(cf. \cite{herz}, Prop. 2.3.1, p.790). $\Gamma(L)$ acts
sharply 3-transitively on $\textup{PG}(1,p^i)$ and therefore (regarded as a
subset of $\Gamma(D)$) on $B$, too. Hence, the order of the stabiliser
of $B$ is $|\Gamma(L)|=p^i(p^{2i}-1)$.\\
In $\Sigma(L,D)$, three mutually nonparallel points are incident with
exactly one chain (cf. \cite{benz1}, Theorem 1.1, p. 95), and the
chains are precisely the blocks of our divisible design.
By Prop. \ref{spera}, we get a $3-(s,k,\lambda_3)$-DD with $s=q$, $k:=|B|
=p^i+1$ and $\lambda_3=1$.\\

\footnotetext{The number refers to the corresponding part of the proof
  below.}

This is also a 2-$(q,p^i+1,\frac{q(q-1)}{p^i-1})$-divisible design.

By removing a set $M\subset B$, we define $B':=B\setminus
M$ as our base block to construct the DD's of the cases (ii)-(v).
The block $B'$ should contain at least three points, therefore, $p^i$
has to be big enough.\\
The stabiliser of $B'$ in $\Gamma(D)$ has to be in $\Gamma(L)$ (see
above) and thus to fix $B$ and hence also $M$ setwise. Therefore,
\[\Gamma(D)_{B'}=\Gamma(L)_M.\]
Notice that, vice versa, the stabiliser of $M$ in $\Gamma(D)$ is equal
to the stabiliser of $B'$ in $\Gamma(D)$ iff $M$ consists of at least
three elements.

\item[(ii)]

Let $B':=B\setminus\{\mathcal{R}(1,0)\}$. This is an $L$-chain minus
one point.
$\Gamma(L)$ acts sharply 3-transitively on $B$, hence
$|\Gamma(L)|=|\Gamma(L)_{\mathcal{R}(1,0)}|(p^i+1)$, and therefore
$|\Gamma(L)_{B'}|=p^i(p^i-1)$.\\
By Proposition \ref{spera}, we get a $3-(s,k,\lambda_3)$-divisible
design with  $s=q$,
\[k:=|B'|=p^i \mbox{ and } \lambda_3=\frac{|\Gamma|}{|\Gamma_{B'}|}
\left(\begin{array}{c} p^i\\3
  \end{array}\right)/\left[\left(\begin{array}{c}\frac{q^2+q}{q}\\ 3
    \end{array}\right)q^3\right]= p^i-2.\]

This is also a 2-$(q,p^i,q(q-1))$-divisible design.

\item[(iii)]

Let $B':=B\setminus\{\mathcal{R}(1,0), \mathcal{R}(0,1)\}$.\\
From the sharp 3-transitivity of $\Gamma(L)$ and the number of
possible permutations of the elements of $M$, we know
$|\Gamma(L)_{\{\mathcal{R}(1,0), \mathcal{R}(0,1)\}}|=2(p^i-1)$. Therefore, by
Proposition \ref{spera}, we get a $3-(s,k,\lambda_3)$-divisible design
with  $s=q$, $k:=|B'|=p^i-1$ and $\lambda_3=\frac{1}{2}(p^i-2)(p^i-3)$.

This is also a 2-$(q,p^i-1,\frac{1}{2}(p^i-2)q(q-1))$-divisible design.

\item[(iv)]

Let $B':=B\setminus\{\mathcal{R}(1,0), \mathcal{R}(0,1), \mathcal{R}(1,1)\}$.\\
Similar to (iii), we can conclude $|\Gamma(L)_{\{\mathcal{R}(1,0),
  \mathcal{R}(0,1), \mathcal{R}(1,1)\}}|=6$ from the sharp
3-transitivity of $\Gamma(L)$ and the fact
that there exist $3!=6$ possible permutations of the elements of
$M$. Now, we get a $3-(s,k,\lambda_3)$-divisible design with  $s=q$,
$k:=|B'|=p^i-2$ and $\lambda_3=\frac{1}{6}(p^i-2)(p^i-3)(p^i-4)$.

This is also a 2-$(q,p^i-2,\frac{1}{6}(p^i-2)(p^i-3)q(q-1))$-divisible
design.

\item[(v)]

Let $B':=B\setminus\{\mathcal{R}(1,0), \mathcal{R}(0,1),
\mathcal{R}(1,1), \mathcal{R}(x,1)\}$, with $x \in L\setminus\{0,1\}$.\\
In analogy to the cases above, the stabiliser of $B'$ in $\Gamma(L)$
corresponds to the stabiliser of $M$ in $\Gamma(L)$. Two 4-tuples of
points are projectively equivalent iff their cross-ratios are
equal. The four points of $M$ allow 24 permutations, but they
determine only the following six cross-ratios:
\[x,\; 1/x,\; 1-x,\; 1/(1-x),\; (x-1)/x,\; x/(x-1) \qquad\hspace*{\fill}(*)\]
In any case, the cross-ratio of the four points is invariant under a
projective group of order 4 isomorphic to
$\mathbb{Z}_2\times\mathbb{Z}_2$ (\cite{hi1}, p.119/120).\\

\begin{itemize}
\item[(a)]
If all six values of $(*)$ are different, then the stabiliser of $B'$
consists only of these four elements since $24=6 \cdot 4$. This case
occurs if the four points form neither a harmonic nor an
equianharmonic quadruple.\\
If $p$ is even, no harmonic quadruple exists, so $x$ can be any
element of $\mathrm{GF}(p^i)\setminus\{0,1\}$ which is not a solution of
$x^2-x+1=0$. We can choose a suitable point $\mathcal{R}(x,1)$ since there
exist at most two solutions of this equation and since $p^i>5$ is
assumed.\\
To get such a point if $p$ is odd, we have to assume that $x$ is
neither an element of $\{0,1,-1,1/2,2\}$ nor a solution of
$x^2-x+1=0$; hence, $p^i>7$ is sufficient.

\item[(b)]
If at least two of the values of $(*)$ are equal, the four points form
a harmonic quadruple if the values of the cross-ratios are $\{-1,1/2,2\}$ or
an equianharmonic one if $x=1/(1-x)$ (or equivalently $x=(x-1)/x$) or both,
which is called superharmonic by Hirschfeld (\cite{hi1})
and which occurs if and
only if $p=3$. In this case, the stabiliser of $B'$ is the symmetric
group $S_4$ of order 24.\\
A harmonic quadruple where $p> 3$ is stabilized by the dihedral group
$D_4$ of order 8. Equianharmonic quadruples exist precisely when $p^i
\equiv 1 \mathrm{mod} 3$ and their stabiliser is the alternating group $A_4$ of
order 12 (\cite{hi1}, p.121).\\
\end{itemize}

\noindent Let $B'$ be a 4-subset of $B$, then we get the same groups as
stabilisers as above. By Proposition \ref{spera}, we get the divisible
designs of Theorem 3.1.1, part (v).\hspace*{\fill}$\Box$
\end{itemize}
\end{paragraph}

\begin{rem}
A.G. Spera constructed divisible designs from a finite local $K$-
algebra $A$ with $K=\mathrm{GF}(q)$ and $J$ its Jacobson radical (with
$|A|=q^n, |J|=q^j,\;n,j \in \Nset$). In the special case $K\cong A/J$,
where $A$ is a Laguerre algebra,  he obtained a transversal
$3-(q^j,q+1,1)$-\textup{DD} as in case (i) (cf. \cite{sp2}).
\end{rem}

\section{Comparing both constructions}

In the introduction, we already mentioned a connection between our
construction and that of Cerroni and Schulz \cite{cesch1}. Since this
connection is not obvious, we mention another representation of the
Laguerre geometry $\Sigma(\mathrm{GF}(q),\mathbb{D}(\mathrm{GF}(q)))$.
After that, we give a short description of Cerroni and Schulz's
construction which will show the duality.\\

Similar to Blaschke's Cylinder-Model (\cite{blasch}, \cite{benz1})
in the real 3-space, it is
possible to embed $\Sigma(K, D)$ in a 3-dimensional projective
space $\Psi$ (cf. \cite{hot}).
By using the more general case showed by Hotje \cite{hot},
we can identify the elements of $\mathbb{P}(D)$ with the elements of a
quadratic cone $\mathcal{O}$, except its vertex $E$, in the
3-dimensional projective space $\Psi$. Similar to Blaschke's
Cylinder-Model, two points are parallel iff they lie on the same generator
\footnote{A generator of $\mathcal{O}$ is a line which is
  completely contained  in $\mathcal{O}$.} of $\mathcal{O}$ \cite{hot}.
Consider a plane in $\Psi$, whose intersection with $\mathcal{O}$ is
exactly the point $E$, then all points of the cylinder
$Z:=\mathcal{O}\setminus E$ are affine points of the 3-dimensional
affine space whose ideal plane is the plane considered above. Such a
plane exists since no finite field is quadratically closed.
$Z$ consists of $q+1$ lines (generators) each containing $q$ points
and intersecting the ideal plane in $E$. Each line contains
precisely one parallel class of points.\\

Now, keeping this in mind, we turn to Cerroni and Schulz's construction
\cite{cesch1}.
Consider a conic $O$ in the ideal plane $E'$ of the 3-dimensional
affine space  AG$(3,q)$. There is a unique tangent at each point of $O$
which determines precisely one  parallel class of affine
planes. Planes of the same parallel class all intersect $E'$ in the
appropriate  tangent. In this way, one gets $q+1$ parallel classes
each consisting of $q$ planes which are the points of the constructed
divisible designs.
By dualising $\Psi$, we obtain the situation of this construction in
which the plane $E'$ is dual to the point $E$ and the planes of one
parallel class correspond to the points of one generator of $Z$,
respectively. The series of divisible designs of Theorem \ref{meinDD},
part (i) possess the same parameters as the
series constructed by Cerroni and Schulz (cf. Theorem
\ref{schulz/cerroni}). They seem to be mutually dual, whereas the other
series of divisible designs of Theorem \ref{meinDD}, arising from different
base blocks, seem to be new.\\

%\bibliographystyle{plain}
%\bibliography{Designs}

Authors' addresses: Sabine Giese and Ralph-Hardo Schulz, 2.~Mathematisches
Institut, FB Mathematik und Informatik, Freie Universtit\"at Berlin, Arnimallee
3, D-14195 Berlin, Germany\\ \emph{giese@mi.fu-berlin.de,
schulz@mi.fu-berlin.de}
\par
Hans Havlicek, Abteilung f\"ur Lineare Algebra und Geometrie, Technische
Universit\"at Wien, Wiedner Hauptstrasse 8-10/1133, A-1040 Wien, Austria\\
\emph{havlicek@geometrie.tuwien.ac.at}

\end{document}